\documentclass{article}
\usepackage[english]{babel}
\usepackage{amsmath,amssymb,bbm,latexsym,tikz-cd}

\newcommand{\tmem}[1]{{\em #1\/}}
\newcommand{\tmop}[1]{\ensuremath{\operatorname{#1}}}
\newcommand{\tmstrong}[1]{\textbf{#1}}
\newcommand{\tmtextbf}[1]{{\bfseries{#1}}}
\newcommand{\tmtextit}[1]{{\itshape{#1}}}
\newenvironment{proof}{\noindent\textbf{Proof\ }}{\hspace*{\fill}$\Box$\medskip}
\newtheorem{theorem}{Theorem}

\begin{document}

\title{Local Coefficients Revisited}

\author{Slawomir Kwasik\footnote{The first author is supported by the Simons Foundation grant 281810.
2010 Mathematics Subject Classification 55P65(primary) 57P10(Secondary).}  and Fang Sun}

\maketitle

\begin{abstract}
Two simple "simplicial approximation" tricks are invoked to prove basic results involving (co)-homology with local coefficients.
\end{abstract}
\section{Introduction}

Homology and cohomology with local coefficients have applications in a variety
of topics including Obstruction Theory, Spectral Sequences, Generalized
Poincare Duality and more. These homology and cohomology theories possess many
properties analogous to those of homology and cohomology with constant
coefficients. Yet some of the corresponding properties are missing. In
particular, there is no Universal Coefficient Theorem linking homology with
local coefficients with cohomology (there is a version in [7], p. 283, though
its application is limited). This (among others) poses extra difficulties in
proving many familiar (and useful) properties of these theories. In fact more than often the proof
of a theorem, carried verbatim from the theory with constant coefficients, is
fundamentally distinct in the case of local coefficients.

The purpose of this paper is to establish three basic properties of
(co)-homology with local coefficients.

\

{\tmstrong{{\underline{The First Property}}}}: A weak homotopy equivalence
induces an isomorphism on cohomology with local coefficients.

\

{\tmstrong{{\underline{The Second Property}}}}: For a CW-pair the long exact
sequences for singular and cellular cohomology group with local coefficients are
equivalent.

\

Despite being very basic, no written account of these properties could be
located. Quite likely this was due to the suspicion that the technical difficulties
in rigorous proofs presented a price too high to pay for the final
product. Our proofs of these properties are based on a simple ``simplicial
approximation'' trick which avoids most of the technicalities.

\

{\tmstrong{{\underline{The Third Property}}}}: The Poincare Duality with
local coefficients for closed, topological manifolds.\\

In this case the written (and complete) accounts of this property do exist. In
fact, we are aware of two such accounts (cf. [8], [9]). Both of these accounts are surprisingly lengthy and technically quite demanding to say the least.

Our argument once again relies on a simple ``simplicial approximation'' trick
of some what different nature than the one mentioned earlier. This leads to a conceptually satisfying and direct proof.

More comments of historical, motivational and mathematical nature are
contained in corresponding sections dealing with the proofs.

\

\section{Local Coefficient Systems}

For more detail on this topic, the reader is referred to [11], Chapter
VI, Section 1-4.

Let $G$ be a bundle of Abelian groups (local coefficient system) on a
topological space $X$. Recall that the{\tmem{ singular chain complex of $X$
with coefficient in $G$}} is defined as
\begin{eqnarray}
  & S_k (X ; G) = \underset{\sigma : \Delta^k \rightarrow X}{\oplus} G
  (\sigma (e_0)) &  \nonumber\\
  & \partial (g \sigma) = G (\bar{\sigma}) (g) \sigma_0 + \underset{i =
  1}{\overset{k}{\Sigma}} (- 1)^i g \sigma_i & 
\end{eqnarray}
where $\Delta^k = \langle e_0, e_1, \cdots, e_k \rangle$ is the standard
$k$-simplex, $g \sigma \in S_k (X ; G)$ is the element that is $g$ on the $G
(\sigma (e_0))$ factor and $0$ otherwise, $\bar{\sigma}$ is $\sigma$ composed
with the straight path from $e_1$ to $e_0$ in $\Delta^k$, $\sigma_i$ is the
restriction of $\sigma$ to the $i$-th face.

Note that the notations we are using are slightly different from the one used
in Whitehead' book.

It can be shown that $\{ S_{\ast} (X ; G), \partial \}$ is a chain complex,
and its homology is called the {\tmem{singular homology}} {\tmem{of}} $X$
{\tmem{with coefficient in}} $G$, denoted as $H_{\ast} (X ; G)$.

In a similar way, one could define the {\tmem{singular cochain complex}}
$S^{\ast} (X ; G)$
\begin{eqnarray}
  & S^k (X ; G) = \underset{\sigma : \Delta^k \rightarrow X}{\Pi} G (\sigma
  (e_0)) &  \nonumber\\
  & (- 1)^k (\delta c) (\sigma) = G (\bar{\sigma})^{- 1} c (\sigma_0) +
  \underset{i = 1}{\overset{n + 1}{\Sigma}} (- 1)^i c (\sigma_i) & 
\end{eqnarray}
where $c \in S^k (X ; G)$.

$\{ S^{\ast} (X ; G), \delta \}$ forms a cochain complex and its homology is
called the {\tmem{singular cohomology}} {\tmem{of}} $X$ {\tmem{with
coefficient in}} $G$, denoted as $H^{\ast} (X ; G)$.

Long exact sequence of pairs, homotopy invariance, excision and additivity
(with respect to disjoint union) are still valid as in the case of constant
coefficients (cf. [11]. VI. 2). As is the equivalence between singular
and cellular homology/cohomology (cf. [11]. VI. 4).

\

\section{Weak Homotopy Equivalence and Cohomology}

In this section, we will provide a proof of the following result:

\begin{theorem}
  Let $f : X \rightarrow Y$ be a weak homotopy equivalence, $G$ be a
  coefficient system on $Y$, then $f_{\ast} : H_{\ast} (X ; f^{\ast} G)
  \rightarrow H_{\ast} (Y ; G), f^{\ast} : H^{\ast} (Y ; G) \rightarrow
  H^{\ast} (X ; f^{\ast} G)$ are isomorphisms, where $f^{\ast} G$ is the
  pull-back of $G$ via $f$. 
\end{theorem}

\

\subsection{Background}

It is a standard result in homotopy theory that weak homotopy equivalences
(continuous maps which induce isomorphisms of all homotopy groups with all
choices of basepoints) induce isomorphism on singular homology and cohomology.
One may ask whether this is true for homology and cohomology with local
coefficients. This result is needed for the following definition in
Obstruction Theory.

Suppose $(K.L)$ is a relative CW-complex, $p : X \rightarrow B$ is a fibration
with $(n - 1)$-connected fiber $F$ and we are given a commutative diagram:\\

\begin{tikzcd}
L \arrow[d, hook] \arrow[r, "f"] & X \arrow[d, "p"]\\
K \arrow[r, "\phi"] & B
\end{tikzcd}\\

The diagram induces an element $\bar{\gamma}^{n + 1} (f) \in H^{n + 1} (K, L
; \phi^{\ast} \pi_n (\mathcal{F}))$, called the {\tmem{primary obstruction}}
{\tmem{to extending}} $f$ ([11], p.298). The name comes from the fact
that $f$ can be extended to a partial lifting $f_{n + 1} : K^{n + 1}
\rightarrow X$ of $\phi$ if and only if $\bar{\gamma}^{n + 1} (f) = 0$.

Some times it is useful to have primary obstruction defined when $(K, L)$ is
replaced by an arbitrary pair $(P, Q)$. In particular, one needs such
definition when defining the Whitney class of a vector bundle over an
arbitrary base (not necessarily homotopic to a CW-complex), or constructing
the Leray-Serre spectral sequences of a fibration over an arbitrary base.

To this end, we can take a CW-approximation $\varphi : (K, L) \rightarrow (P,
Q)$, i.e., a map of pairs such that $(K, L)$ is a CW-pair and $\varphi,
\varphi_{|L}$ are both weak homotopy equivalences. Thus we have a diagram:\\

\begin{tikzcd}
L \arrow[r, "\varphi"] \arrow[d, hook] & Q \arrow[d, hook] \arrow[r, "f"] & X \arrow[d, "p"]\\
K \arrow[r, "\varphi"] & P \arrow[r, "\phi"] & B
\end{tikzcd}\\

The element $\bar{\gamma}^{n + 1} (f \circ \varphi) \in H^{n + 1} (K, L ;
\varphi^{\ast} \phi^{\ast} \pi_n (\mathcal{F}))$ is well-defined. If we know
that
\[ \varphi^{\ast} : H^{n + 1} (P, Q ; \phi^{\ast} \pi_n (\mathcal{F}))
   \longrightarrow H^{n + 1} (K, L ; \varphi^{\ast} \phi^{\ast} \pi_n
   (\mathcal{F})) \]
is an isomorphism (in [11] p.300 this is assumed without any explanation), then one could define $\bar{\gamma}^{n + 1} (f) : =
\varphi^{\ast - 1} \bar{\gamma}^{n + 1} (f \circ \varphi) \in H^{n + 1} (P, Q
; \phi^{\ast} \pi_n (\mathcal{F}))$. An easy argument of naturality shows that
this is independent of the CW-approximation $\varphi$. Hence $\bar{\gamma}^{n
+ 1} (f)$ is well-defined.

\

There are at least two ways to prove that a weak homotopy equivalence induces
isomorphism on homology (with constant coefficients) in the literature. One
approach uses Hurewicz Theorem ([7] 7.5.9, 7.6.25), the other proof
([2], Proposition 4.21) is by a construction that relies heavily on the
finiteness of singular chains. The analogous result for cohomology with
constant coefficients follows from this via the Universal Coefficient Theorem.

\begin{table}[h]
  \begin{tabular}{ccc}
    & Constant coefficients & Local coefficients\\
    Homology & Hurewicz/Construction & Construction\\
    Cohomology & Universal Coefficient Theorem & ?
  \end{tabular}
  \caption{}
\end{table}

When it comes to local coefficients, Hurewicz Theorem is no longer available.
The constructive proof still works for homology with local coefficients. Yet
due to the absence of Universal Coefficient Theorem, the result for cohomology
does not follow automatically. Our proof turns out to be quite different from
those above. As far as we know, no alternative exists in the literature for
cohomology.

\

\subsection{Singular Complex}

We begin with the notion of singular complex, which is central in our proof.

Let $X$ be a topological space. Take the disjoint union of $k$-simplexes, one
for each continuous map $\sigma : \Delta^k \rightarrow X$. Do this for all
integer $k \geqslant 0$. Glue the simplexes according to restriction of maps
to faces. The resulted CW-complex is called the singular complex of $X$,
denoted by $S X$. In fact, $S X$ is a $\Delta$-complex ([2], p.164).

A continuous map $f : X \rightarrow Y$ induces (with obvious definition) $S f
: S X \rightarrow S Y$.

Since the simplexes of $S X$ corresponds to continuous maps $\Delta^k
\rightarrow X$, there is a canonical map $I_X : S X \rightarrow X$ mapping
each simplex via the map defining it. $I_X$ is natural with respect to
continuous map $f : X \rightarrow Y$, i.e., the following diagram commutes:\\

\begin{tikzcd}
X \arrow[r, "f"] & Y\\
S X \arrow[u, "I_X"] \arrow[r, "Sf"] & S Y \arrow[u, "I_Y"]
\end{tikzcd}\\\\

The following result will be useful for our purpose:\\

{\noindent}\tmtextbf{Theorem 2 }\tmtextit{For any topological space X, $I_X$
is a weak homotopy equivalence.}{\hspace*{\fill}}{\medskip}

\begin{proof}
  See [5], Chapter III, Theorem 6.7, or [1] Theorem 16.43 on p.149.
\end{proof}

\

For a coefficient system $G$ on a $\Delta$-complex $K$, there is a version of
simplicial homology/cohomology. The definition of chains/cochains and
boundary/coboundary maps are the same with that of singular ones, except that
the direct sum/product is now over all simplices of $K$. We denote the
simplicial chain/cochain of $K$ with coefficient in $G$ by $\Delta_{\ast} (K ;
G)$ and $\Delta^{\ast} (K ; G)$, and simplicial homology/cohomology by
$H_{\ast}^{\Delta} (X ; G)$ and $H_{\Delta}^{\ast} (X ; G)$.

As in the case of constant coefficients, we still have:

{\noindent}\tmtextbf{Theorem 3 }\tmtextit{The canonical injection
$\Delta_{\ast} (K ; G) \overset{j_{\#}}{\hookrightarrow} S_{\ast} (K ; G)$ and
projection $S^{\ast} (K ; G) \overset{j^{\#}}{\rightarrow} \Delta^{\ast} (K ;
G)$ are chain $I_X^{\ast} : H^{\ast} (X ; G) \rightarrow H_{\Delta}^{\ast} (S
X ; I_X^{\ast} G)$maps that induce isomorphism on homology/cohomology
groups.}{\hspace*{\fill}}{\medskip}

\begin{proof}
  The proof for homology is the same as in the case of constant coefficients
  (cf. [2], 2.27). Although Universal Coefficient Theorem is not
  available, one could easily adapt the above mentioned proof to the case of
  cohomology, with little change.
\end{proof}

Now suppose $G$ is a coefficient system on $X$. Let $I_X^{\ast} G$ be the
pullback of $G$ via $I_X$. We have
\[ S^{\ast} (X ; G) \overset{I_X^{\#}}{\longrightarrow} S^{\ast} (S X ;
   I_X^{\ast} G) \overset{j\#}{\longrightarrow} \Delta^{\ast} (S X ;
   I_X^{\ast} G) \]
It is easy to see that the composition $j^{\#} \circ I_X^{\#}$ identifies $S^k
(X ; G) = \underset{\sigma : \Delta^k \rightarrow X}{\Pi} G (\sigma (e_0)) =
\Delta^{\ast} (S X ; I_X^{\ast} G)$. In particular, $j^{\#} \circ I_X^{\#}$ is
a chain isomorphism and $j^{\ast} \circ I_X^{\ast} : H^{\ast} (X ; G)
\rightarrow H_{\Delta}^{\ast} (S X ; I_X^{\ast} G)$ is also an isomorphism. By
a similar argument one can show that $I_{X \ast} \circ j_{\ast} :
H_{\ast}^{\Delta} (S X ; I_X^{\ast} G) \rightarrow H_{\ast} (X ; G)$ is an
isomorphism. Combined with Theorem 3 we have shown:\\

{\noindent}\tmtextbf{Theorem 4 }\tmtextit{The map $I_X$ induces isomorphisms
on homology and cohomology. To be more precise, for any coefficient system $G$
on X, $I_{X \ast} : H_{\ast} (S X ; I_X^{\ast} G) \rightarrow H_{\ast} (X ;
G)$ and $I_X^{\ast} : H^{\ast} (X ; G) \rightarrow H^{\ast} (S X ; I_X^{\ast}
G)$ are isomorphisms.}{\hspace*{\fill}}{\medskip}

\

\subsection{Proof of Theorem 1}

\begin{proof}
  Let $f : X \rightarrow Y$ be a weak homotopy equivalence. As noted above
  there's a commutative diagram

\begin{tikzcd}
X \arrow[r, "f"] & Y\\
S X \arrow[u, "I_X"] \arrow[r, "Sf"] & S Y \arrow[u, "I_Y"]
\end{tikzcd}

\begin{flushleft}
in which $I_X, I_Y$ are weak homotopy equivalences by Theorem 2.
\end{flushleft}
  Commutativity implies that $S f$ is also a weak homotopy equivalence. Since
  $S X, S Y$ are CW-complexes, Whitehead Theorem (cf. [2]) implies that $S f$ is a
  homotopy equivalence, hence induce isomorphism on homology/cohomology with
  local coefficients.
  
  Now apply Theorem 4 to the induced (commutative) diagram on
  homology/cohomology we have the desired result.
\end{proof}

\section{Identifying Singular and Cellular Long Exact Sequences}

\subsection{Background}

Let $(K, L)$ be a CW-pair, $G$ be a coefficient system on $K$. There is long
exact sequence
\begin{equation}
  \longrightarrow H_n (L ; G)
  \longrightarrow H_n (K ; G) \longrightarrow H_n (K, L ; G) \longrightarrow
\end{equation}
and
\begin{equation}
  \longrightarrow H^n (K, L ; G) \longrightarrow H^n (K ; G) \longrightarrow
  H^n (L ; G) \longrightarrow
\end{equation}
There are also naturally defined short exact sequences of cellular
chain/cochain complexes
\[ 0 \longrightarrow \Gamma_{\ast} (L ; G) \longrightarrow \Gamma_{\ast} (K ;
   G) \longrightarrow \Gamma_{\ast} (K, L ; G) \longrightarrow 0 \]
and
\[ 0 \longrightarrow \Gamma^{\ast} (K, L ; G) \longrightarrow \Gamma^{\ast} (K
   ; G) \longrightarrow \Gamma^{\ast} (L ; G) \longrightarrow 0 \]
which induce long exact sequences
\begin{equation}
 \longrightarrow H_n (\Gamma_{\ast} (L)
  ; G) \longrightarrow H_n (\Gamma_{\ast} (K) ; G) \longrightarrow H_n
  (\Gamma_{\ast} (K, L) ; G) \longrightarrow
\end{equation}
and
\begin{equation}
 \longrightarrow H^n (\Gamma^{\ast} (K, L) ; G) \longrightarrow H^n (\Gamma^{\ast} (K) ; G)
  \longrightarrow H^n (\Gamma^{\ast} (L) ; G) \longrightarrow
\end{equation}

The groups in (3) and (5) (resp. (4) and (6)) are term-wise isomorphic. It is
natural to ask whether the the long exact sequences (viewed as chain complexes) are chain isomorphic. This
is used in the proof of [11] Theorem VI.6.9 (again this is used without any justification). It is natural to expect that this problem can be solved
by a diagram chasing, since cellular homology/cohomology are themselves
defined by certain diagrams. Yet as far as we know, no such proof has been
given. In fact the only relevant result in the literature is given in Schubert' book (cf.
[6], p. 303). Schubert constructed an intermediate between the singular and
cellular chain complex, called normal chain complex and used it to show (3)
and (5) are chain isomorphic. The construction (again!) depends heavily on the
finiteness of singular chains, thus fails to prove the result for cohomology
with local coefficients (though for constant coefficients one could still use
Universal Coefficient Theorem to dualize everything).

Our goal is to prove:\\

{\noindent}\tmtextbf{Theorem 5 }\tmtextit{The long exact sequences (3) and (5)
(resp. (4) and (6))are chain isomorphic.}{\hspace*{\fill}}{\medskip}

We shall prove the result for cohomology, the proof for homology is analogous.

\

\subsection{Proof of Theorem 5}

\begin{proof}
  Since the identification of singular and cellular homology/cohomology groups
  are natural with respect to cellular maps (whether the coefficient is
  constant or local), it suffice to check the commutativity of the diagram
  (coefficients omitted):\\
  
  \begin{tikzcd}
H^n (\Gamma^* (L)) \arrow[r, "\delta"] \arrow[d, leftrightarrow] & H^{n+1} (\Gamma^* (K,L))\arrow[d, leftrightarrow]\\
H^n (L) \arrow[r, "\delta"] & H^{n+1} (K,L)
\end{tikzcd}\\

  Note that $S L$ can be identified canonically with a subspace of $S K$ and $I_K : S K
  \rightarrow K$ restricts to $I_L$ on $S L$. Thus we have a weak homotopy
  equivalence $I_K : (S K, S L) \rightarrow (K, L)$. Since both the domain and
  codomain are CW-complexes, $I_K$ is actually a homotopy equivalence.
  Homotope $I_K$ to a cellular map $J_K$ and consider:\\

  \begin{tikzcd}[row sep=scriptsize, column sep=scriptsize]
& H^n (L)\arrow[dl] \arrow[rr] \arrow[dd] & & H^{n+1} (K,L) \arrow[dl] \arrow[dd] \\
H^n (S L) \arrow[rr, crossing over] \arrow[dd] & & H^{n+1} (S K, S L)  \\
& H^n (\Gamma^* (L)) \arrow[dl] \arrow[rr] & & H^{n+1} (\Gamma^* (K,L)) \arrow[dl] \\
H^n (\Gamma^* (S L)) \arrow[rr] & & H^{n+1} (\Gamma^* (S K,S L)) \arrow[from=uu, crossing over]\\
\end{tikzcd}
  in which vertical arrows are isomorphism between singular and cellular
  cohomology groups, horizontal arrows are boundary maps in the corresponding
  long exact sequence and arrows going down left are induced by $J_K$.
  Coefficients are obvious and omitted.
  
  The rectangles on top of the above diagram commute since a map ( in this
  case $J_K$) induces a chain map between singular long exact sequences.
  Similarly, the cellular map $J_K$ induce chain map between cellular long
  exact sequences, hence the commutativity of the bottom rectangle. The
  rectangles on the left and right commute thanks to naturality of the
  identification of singular cohomology with cellular cohomology under
  cellular maps.
  
  All down left arrows are isomorphisms since $J_K$ is a (cellular) homotopy
  equivalence. In particular the map $J_K^{\ast} : H^{n + 1} (\Gamma^{\ast}
  (K, L)) \rightarrow H^{n + 1} (\Gamma^{\ast} (S K, S L))$ is an isomorphism.
  
  We intend to prove the commutativity of the rectangle in the back. As
  indicated by the above argument, it suffice to show that for the front
  rectangle. In other words, we have reduced the problem to the case where
  $(K, L)$ is a pair of $\Delta$-complexes. We shall assume this from now.
  
  For a $\Delta$-complex pair $(K, L)$ and a coefficient system $G$ on $K$,
  there is an isomorphism $\Phi : \Delta^n (K, L ; G) \rightarrow \Gamma^n (K,
  L ; G)$ defined by the identification
  \[ \Delta^n (K, L ; G) \longleftrightarrow \underset{\sigma}{\Pi} G (\sigma
     (e_0)) \longleftrightarrow \underset{\sigma}{\Pi} H^n (\Delta^n, \partial
     \Delta^n ; \sigma^{\ast} G) \longleftrightarrow \Gamma^n (K, L ; G) \]
  where the direct products are over all $n$-simplexes of $K - L$. It is easy
  to check (by a diagram chasing) that $\Phi$ commute with boundary maps of
  the two chain complexes and hence is a chain isomorphism.
  
  $j^{\#}$ and $\Phi$ induces the following commutative diagram joining the
  singular, simplicial, and cellular short exact sequences of $(K, L)$
  (coefficients omitted)\\

  \begin{tikzcd}
0 \arrow[r] & S^*(K,L) \arrow[r] \arrow[d] & S^*(K)\arrow[r] \arrow[d] & S^*(L) \arrow[r] \arrow[d] & 0\\
0 \arrow[r] & \Delta^*(K,L) \arrow[r] \arrow[d] & \Delta^*(K)\arrow[r] \arrow[d] & \Delta^*(L) \arrow[r] \arrow[d] & 0\\ 
0 \arrow[r] & \Gamma^*(K,L) \arrow[r] & \Gamma^*(K)\arrow[r] & \Gamma^*(L) \arrow[r] & 0\\
\end{tikzcd}
\begin{flushleft}
which induces a commutative diagram for the boundary homomorphisms in the
  corresponding long exact sequences
  \end{flushleft}

  \begin{tikzcd}
H^n (L) \arrow[r, "\delta"] \arrow[d, "j^*"] & H^{n+1} (K,L)\arrow[d,"j^*"]\\
H^n_\Delta (L) \arrow[r, "\delta"] \arrow[d, "\Phi^*"] & H^{n+1}_\Delta (K,L) \arrow[d, "\Phi^*"]\\
H^n (\Gamma^* (L)) \arrow[r, "\delta"] & H^{n+1} (\Gamma^* (K,L))
\end{tikzcd}\\

  The theorem will follow from the lemma below.
\end{proof}

{\noindent}\tmtextbf{Lemma 6 }\tmtextit{The isomorphism $\Phi^{\ast} \circ
j^{\ast}$ is exactly the canonical identification between cellular and
singular cohomology.}{\hspace*{\fill}}{\medskip}

\begin{proof}
  We shall prove this for the absolute case (i.e. $L = \varnothing$). The
  proof of the relative case is similar. The coefficient system $G$ will be
  omitted unless necessary.
  
  Consider the following diagram:\\
  
\begin{tikzcd}
H^n (K) \arrow[rd, "j^*"] \arrow[d, "\iota"]\\
H^n (K^n) &  H^n_\Delta (K) \arrow[ld, "\Phi^*"] \\
H^n (\Gamma^*(K)) \arrow[u, "\kappa"]
\end{tikzcd}\\

  where the $\iota$ is induced by the inclusion $K^n \hookrightarrow K$ and
  $\kappa$ is induced by the homomorphism (induced by identity) $\Gamma^n (K)
  = H^n (K^n, K^{n - 1}) \rightarrow H^n (K^n)$.
  
  For any $[b] \in H^n (K)$, where $b \in S^n (K)$. Let $[\alpha] \in H^n
  (\Gamma (K))$ be the element corresponding to $[b]$ via the canonical
  identification between singular and cellular cohomology. Then $\kappa
  ([\alpha]) = \iota ([b])$, hence $\alpha = [a]$ for some $a \in S^n (K^n,
  K^{n - 1}) \subset S^n (K^n)$ such that $a - b_{|K^n} = \delta c$ for some
  $c \in S^{n - 1} (K^n)$.
  
  It suffice then to show that $\Phi^{\ast - 1} ([\alpha]) = j^{\ast} ([b])$.\\

  We know that
   \begin{align}
    & j^{\#} (a) - j^{\#} (b) & \nonumber \\
    = & j^{\#} (a - b_{|K^n}) & \nonumber\\
    = & j^{\#} (\delta c) & \nonumber\\
    = & \delta j^{\#} (c) \in \Delta^n (K^n) = \Delta^n (K) & \nonumber
  \end{align}
  In other words, $j^{\#} (a), j^{\#} (b)$ are cohomologous.\\

  Also, one can check that
  \begin{align}
    & j^{\#} (a) & \nonumber\\
    = & \Phi^{\#- 1} ([a]) & \nonumber\\
    = & \Phi^{\#- 1} (\alpha) \in \Delta^n (K) = \underset{\sigma}{\Pi} G
    (\sigma (e_0)) \nonumber
  \end{align}
  by looking at their value on each $n$-simplex $\sigma$.
  
  Thus $\Phi^{\ast - 1} ([\alpha]) = [\Phi^{\#- 1} (\alpha)] = [j^{\#} (a)] =
  [j^{\#} (b)] = j^{\ast} ([b])$.
\end{proof}

\

\

\section{Simplicial Approximation and Poincare Duality}

\

We now turn to another type of simplicial approximation. Let $\mathcal{M}$ be
a $n$-dimensional topological manifold, $R$ be a principal ideal domain and
$G$ be a bundle of right $R$-modules.

The Poincare Duality Theorem (with local coefficients) states that
\[ H^i_c (\mathcal{M}; G) \overset{\frown \mu_{\mathcal{M}}}{\longrightarrow}
   H_{n - i} (G \otimes_R \mathcal{M}_R), 0 \leqslant i \leqslant n \]
where $H^{\ast}_c$ stands for singular cohomology with compact support,
$\mu_{\mathcal{M}}$ is the (generalized) fundamental class and $\mathcal{M}_R$
is the orientation bundle of $\mathcal{M}$ with coefficient in $R$.

For relevant definitions and proof of the theorem, see [8] or [9].

There is a version of this duality for compact triangulated manifolds with or without boundary (see
[4], or [10] Theorem 2.1 p.23), which dates back much earlier i.e., the original
proof given by S. Lefschetz. This proof is short and purely geometric (it
uses the dual decomposition of the corresponding simplicial complex). Thus it
would be nice if one could reduce the general case to the case of triangulated manifolds.
This is when simplicial approximation comes into the picture.

Assume, for simplicity, that $\mathcal{M}$ is closed and orientable.\\

{\noindent}\tmtextbf{Theorem 7 }\tmtextit{Let $G$ be a coefficient system on
$\mathcal{M}$, then $H^i (\mathcal{M}; G) \overset{\frown
\mu_{\mathcal{M}}}{\longrightarrow} H_{n - i} (G)$ is an isomophrism for all
$0 \leqslant i \leqslant n$.}{\hspace*{\fill}}{\medskip}

\begin{proof}
  By [3], there is a $k$-disk bundle $p : E \rightarrow
  \mathcal{M}$ such that $E$ admits a triangulation. Also, $E$ is embedded in
  $\mathbbm{R}^{n + k}$. In particular, $E$ is orientable as a manifold with
  boundary. There is a Poincare Duality
  \[ H^{i + k} (E, \partial E ; p^{\ast} G) \overset{\frown
     \mu_E}{\longrightarrow} H_{n - i} (E ; p^{\ast} G) \]
  where $\mu_E$ is the fundamental class of $E$ and $p^{\ast} G$ is the
  pull-back bundle. Since $E$ is triangulable, this is an isomorphism.
  
  Next we prove that $E$ is orientable as a disk bundle, i.e., there exist $U
  \in H^k (E, \partial E) = H^k (E, \partial E ; \mathbbm{Z})$ that restrict
  to a generator $U_{|x} \in H^k (E_x, \partial E_x)$ for every $x \in
  \mathcal{M}$. Here $(E_x, \partial E_x)$ stands for the fiber over $x$.
  
  The composition
  \[ H^k (E, \partial E) \overset{\frown \mu_E}{\longrightarrow} H_n (E)
     \overset{p_{\ast}}{\longrightarrow} H_n (\mathcal{M}) \]
  is an isomorphism since $p$ is a homotopy equivalence. Define $U \in H^k (E,
  \partial E)$ by $U \frown \mu_E = p_{\ast}^{- 1} (\mu_{\mathcal{M}})$.
  
  Note that $\mathcal{M}$ embeds in $E$ by the zero-section. For $x \in
  \mathcal{M}$, consider the following diagram, where vertical maps are
  induced by inclusions\\
  
  \begin{tikzcd}
H^k(E,\partial E) \otimes H_{n+k}(E,\partial E) \arrow[d, "i_1", shift right=5ex] \arrow[r, "\frown"] & H^n (E) \arrow[r, "p_*"] \arrow[d, "i_3"] & H_n(\mathcal{M}) \arrow[d, "i_4"]\\
H^k(E|\mathcal{M}) \otimes H_{n+k}(E|(x,0)) \arrow[u, "i_2", shift right=5ex] \arrow[r, "\frown"] & H^n (E|E_x) \arrow[r, "p_*"] & H_n(\mathcal{M}|x)
\end{tikzcd}\\

  Obviously the rightmost box is commutative, and $i_1$ is an
  isomorphism. Let $U' = i_1^{- 1} (U) \in H^k (E|\mathcal{M})$. Define
  $\mu_{(x, 0)} = i_2 (\mu_E) \in H_{n + k} (E| (x, 0))$. Then by naturality of
  cap products, $U' \frown \mu_{(x, 0)} = i_3 p^{- 1}_{\ast}
  \mu_{\mathcal{M}}$. By commutativity and definition of the fundamental class
  $\mu_{\mathcal{M}}$, $p_{\ast} i_{3 \ast} p^{- 1}_{\ast} \mu_{\mathcal{M}} =
  i_4 \mu_{\mathcal{M}}$ is a generator of $H_n (\mathcal{M}|x)$. Since $p :
  (E|E_x) \rightarrow (\mathcal{M}|x)$ is a homotopy equivalence, $i_{3 \ast}
  p^{- 1}_{\ast} \mu_{\mathcal{M}}$ is also a generator.
  
  Now choose a open disk neighborhood $W$ of $x$ in $\mathcal{M}$ upon which
  $E$ admits a local trivialization $\Phi : p^{- 1} (W) \rightarrow W \times
  D^k$. Let $\Phi_x : E_x \rightarrow D^k$ be the restriction of $\Phi$.
  Consider the following diagram\\

 \begin{tikzcd}
H^k(E|\mathcal{M}) \otimes H_{n+k}(E|(x,0))  \arrow[d,"\Phi_1", shift right=5ex] \arrow[r, "\frown"] & H^n (E|E_x)\\
H^k(W\times D^k|W\times 0) \otimes H_{n+k}(W\times D^k|(x,0)) \arrow[u,"\Phi_2", shift right=5ex] \arrow[r, "\frown"] & H^n (W\times D^k|E_x) \arrow[u,"\Phi_3"]
\end{tikzcd}\\
  where vertical homomorphisms are induced by $\Phi$.
  
  By excision, $\Phi_2, \Phi_3$ are isomorphisms. Hence by naturality $\Phi_1
  (U') \frown \Phi_2^{- 1} (\mu_{(x, 0)}) = \Phi_3^{} i_{3 \ast} p^{-
  1}_{\ast} \mu_{\mathcal{M}}$, which is a generator of $H_n (W \times D^k
  |E_x) \cong \mathbbm{Z}$. Note that cross product induces an isomorphism
  \[ H^k (W \times D^k |W \times 0) \longleftrightarrow H^k (D^k |0) \otimes
     H^0 (W) \]
  and $\Phi_1 (U')$ corresponds to $\Phi_{x \ast} (U_{|x}) \times 1$. This
  forces $\Phi_{x \ast} (U_{|x})$ and thus $U_{|x}$ to be a generator. This
  proves that $U$ is an orientation for the disk bundle $E$.
  
  Now we have a diagram\\
  
\begin{tikzcd}
H^i (\mathcal{M};G) \arrow[r,"\frown \mu_{\mathcal{M}}"] \arrow[d] & H_{n-i} (\mathcal{M};G)\\
H^{i+k} (E,\partial E; p^*G) \arrow[r, "\frown \mu_{E}"] & H_{n-i} (E; p^*G) \arrow[u, "p_{*}"]
\end{tikzcd}\\
  in which the left vertical map is the Thom isomorphism (cf. [7]
  p.283) sending $\alpha$ to $p^{\ast} (\alpha) \smile U$.
  
  This diagram commutes:
  \begin{align}
    & p_{\ast} (((p^{\ast} \alpha) \smile U) \frown \mu_E) & \nonumber\\
    = & p_{\ast} (p^{\ast} \alpha \frown (U \smile \mu_E) ) & \nonumber\\
    = & \alpha \frown p_{\ast} (U \frown \mu_E) & (\tmop{naturality})\nonumber\\
    = & \alpha \frown \mu_{\mathcal{M}} & (\tmop{definition} \tmop{of} U)\nonumber
  \end{align}
  The top arrow is thus an isomorphism since all others are isomorphisms.
\end{proof}

\begin{center}
{\Large References}
\end{center}

\begin{flushleft}
[1] B. Gary, Homotopy Theory. \textit{An introduction to Algebraic Topology}, Academic Press, new York, San Francisco, London, 1975.

[2] A. Hatcher, \textit{Algebraic Topology}, Cambridge University Press, 2002.

[3] R. Kirby and L. Siebermann, \textit{Foundational Essays on Topological Manifolds}, Smoothings and Triangulations, Princeton Univ. Press, Princeton 1977.

[4] J. A. Lees, \textit{The Surgery Obstruction Groups of C.T.C. Wall}, Adv. Math. 11(1973), 113-156.

[5] A. Lundell and S. Weingram, \textit{The Topology of CW complexes}, Van Nostrand Reinhold, 1963.

[6] H. Schubert, \textit{Topology}, Allan and Bacon, Inc., 1968, Boston.

[7] E. Spanier, \textit{Algebraic Topology}, McGraw-Hill, 1966.

[8] E. Spanier, \textit{Singular homology and cohomology with local coefficients and duality for manifolds}, Pacific J. Math. 160 (1993), no. 1, 165-200.

[9] F. Sun, \textit{An elementary proof for Poincare Duality with local coefficients}, arXiv:1709.00569

[10] C.T.C. Wall, \textit{Surgery on Compact Manifolds}, Academic Press, New york, 1970.

[11] G. Whitehead, \textit{Elements of Homotopy Theory}, Springer-Verlag, 1978.
\end{flushleft}

\

Slawomir Kwasik

Department of Mathematics

Tulane University

New Orleans, LA, 70118

kwasik@tulane.edu

\

Fang Sun

Department of Mathematics

Tulane University

New Orleans, LA, 70118

fsun@tulane.edu

\end{document}